\documentclass[a4paper]{article}

\usepackage{MyLaTeX,Myarticle}

\title{On Tensor Products of Simple Modules for Simple Groups}
\author{David A.~Craven, University of Oxford}
\date{January 2011}
\renewcommand{\Alt}[1]{\mathfrak{A}_{#1}}
\renewcommand{\Sym}[1]{\mathfrak{S}_{#1}}

\begin{document}
\maketitle

\noindent In an attempt to get some information on the multiplicative structure of the Green ring we study algebraic modules for simple groups, and associated groups such as quasisimple and almost-simple groups. We prove that, for almost all groups of Lie type in defining characteristic, the natural module is non-algebraic. For alternating and symmetric groups, we prove that the simple modules in $p$-blocks with defect groups of order $p^2$ are algebraic, for $p\leq 5$. Finally, we analyze nine sporadic groups, finding that all simple modules are algebraic for various primes and sporadic groups.

\section{Introduction}

The tensor structure of the category of finite-dimensional $kG$-modules, where $k$ is a field of characteristic $p$ and $G$ is a finite group, is a structure that remains largely shrouded in mystery. As a first approximation to understanding this structure, Alperin introduced the notion of \emph{algebraic} modules \cite{alperin1976b}; a $kG$-module is \emph{algebraic} if it satisfies a polynomial with integer coefficients, where addition and multiplication are given by the direct sum and tensor product. Trivial source modules are algebraic, as are simple modules for $p$-soluble groups \cite{feit1980}, simple modules in characteristic $2$ for finite groups with abelian Sylow $2$-subgroups \cite{craven2009}, and simple modules for $2$-blocks with Klein-four defect groups in general \cite{cekl2008}.

Since simple modules for $p$-soluble groups are algebraic, it is natural to ask about simple modules for simple groups, or groups that are close to simple, such as quasisimple and almost simple groups. In this article, we consider the `natural module' for groups of Lie type in defining characteristic, simple modules for various low-rank Lie type groups in small, non-defining characteristics, alternating and symmetric groups, and nine of the sporadic groups.

On natural modules, we have the following result.

\begin{thma}Let $q$ be a power of a prime, and let $G$ be a finite group of Lie type. If $G$ is a classical group, suppose that it is one of the groups $\SL_n(q)$ or $\Sp_{2n}(q)$ for $n\geq 3$, $\SU_n(q)$ for $n\geq 6$, $\Omega_{2n}^+(q)$ or $\Omega_{2n}^-(q)$ for $n\geq 4$, or $\Omega_{2n+1}(q)$ for $n\geq 3$. If $G$ is of exceptional type, suppose that it is one of the groups $^3D_4(q)$, $^2G_2(q)$, $^2F_4(q)$, $G_2(q)$, $F_4(q)$, $E_6(q)$, $^2E_6(q)$, $E_7(q)$, and $E_8(q)$. Then the natural module for $G$ is non-algebraic.
\end{thma}

For the classical groups, the natural module is clear, and we will define what we mean by natural module for the exceptional groups in Section \ref{sec:exceptional}, but it is generally the non-trivial simple module of smallest dimension.

For the Lie-type groups not mentioned in this theorem ($\SL_2(q)$, $\Sp_4(q)$, $\SU_n(q)$ for $n\leq 5$ and $\Suz(q)$), it is only known for $\SL_2(q)$ that the natural module is actually algebraic: this is a result of Alperin \cite{alperin1979} for $q$ even, and unpublished work of Kov\'acs \cite{kovacs1981un} for $q$ odd. Since the $q$ odd case remains unpublished, here we will provide a short proof using the theory of tilting modules, from which it is possible to calculate any tensor power of the natural module for $\SL_2(q)$.

\begin{thma}[Alperin, Kov\'acs] \label{thm:SL2} Let $p$ be a prime, let $n\geq 1$ be an integer, and let $k$ be a field of characteristic $p$. If $G=\SL_2(p^n)$, then all simple $kG$-modules are algebraic.
\end{thma}

This theorem has the following obvious corollary, using the fact that the natural module for $\SL_2(q)$, hence for $\GL_2(q)$, is algebraic, and restrictions of algebraic modules are algebraic.

\begin{cora}\label{cor:2dim} Let $p$ be a prime, and let $k$ be a finite field. If $G$ is a finite group, and $M$ is a $2$-dimensional $kG$-module, then $M$ is algebraic.\end{cora}

Turning to non-defining characteristic, our results focus on the primes $2$ and $3$. If $\charac k=2$, then we analyze the groups $\PSL_2(q)$, $\PSL_3(q)$ (for $q\equiv 3\bmod 4$) and $\PSU_3(q)$ (for $q\equiv 1\bmod 4$). In other words, we analyze the cases where the Sylow $2$-subgroups are either dihedral or semidihedral. (In the next two theorems, we take an algebraically closed field, purely so that we are guaranteed the existence of all simple modules that `should' exist.)

\begin{thma} Let $k$ be an algebraically closed field of characteristic $2$.
\begin{enumerate}
\item If $G=\PSL_2(q)$, then all simple $kG$-modules are algebraic if and only if $q\not\equiv 7\bmod 8$. If $q\equiv 7\bmod 8$, then the two $(q-1)/2$-dimensional simple modules are non-algebraic.
\item If $G=\PSL_3(q)$, then all simple $kG$-modules are algebraic if $q\equiv 3\bmod 8$. If $q\equiv 7\bmod 8$, then the two non-trivial simple modules in  the principal $2$-block are non-algebraic.
\item If $G=\PSU_3(q)$, then all simple $kG$-modules are algebraic if $q\equiv 1\bmod 4$.
\end{enumerate}
\end{thma}

\begin{thma}\label{thm:c3c3simplegroup} Let $k$ be an algebraically closed field of characteristic $3$, and let $G$ be a finite simple group such that $C_3\times C_3$ is a Sylow $3$-subgroup of $G$.
\begin{enumerate}
\item If $G$ is a finite group of Lie type, an alternating group, or one of the sporadic groups $M_{22}$ or $HS$, then all simple $kG$-modules in the principal $3$-block are algebraic.
\item If $G$ is one of the groups $M_{11}$ or $M_{23}$, then there are non-algebraic simple modules in the principal $3$-block.
\end{enumerate}
\end{thma}

Notice that, if $G$ is a finite group of Lie type with elementary abelian Sylow $3$-subgroups, and $k$ is a field of characteristic $3$, then all simple $kG$-modules in the principal block are algebraic. In characteristic $5$, however, this does not occur, thanks to the following proposition.

\begin{propa}\label{prop:secondrow} Let $k$ be a field of characteristic $p>2$, let $G$ be a finite group with non-cyclic Sylow $p$-subgroups, and let $M$ be an indecomposable $kG$-module with $p\nmid\dim M$. If $M$ lies on the second row of its component of the Auslander--Reiten quiver, then $M$ is non-algebraic.
\end{propa}

For $F_4(2)$ in characteristic $5$, the simple module $M$ of dimension 875823 lies on the second row of its component of the Auslander--Reiten quiver \cite{kawatamichleruno2001}, and the Sylow $5$-subgroups of $F_4(2)$ are non-cyclic, so $M$ is non-algebraic.

\medskip

One can use Theorem \ref{thm:c3c3simplegroup}, together with a standard reduction, to determine exactly which finite groups with Sylow $3$-subgroups $C_3\times C_3$ have algebraic simple modules in the principal $3$-block.

\begin{cora}\label{cor:c3c3group} Let $G$ be a finite group with Sylow $3$-subgroups of order at most $9$, and let $k$ be a field of characteristic $3$. Write $B$ for the principal $3$-block of $kG$. All simple $B$-modules are algebraic if and only if neither $M_{11}$ nor $M_{23}$ is a composition factor of $G$.
\end{cora}

For symmetric and alternating groups, we study simple modules in blocks with defect group $C_p\times C_p$, and prove the following result.

\begin{thma} Let $p$ be one of $2$, $3$ and $5$, and let $k$ be a field of characteristic $p$. Let $G$ be either an alternating group $\Alt{n}$ or a symmetric group $\Sym{n}$ for some $n$. If $B$ is a block of $kG$ with defect group $C_p\times C_p$, then all simple $B$-modules are algebraic.
\end{thma}

The reason that this result cannot be extended further is that it is not possible to prove whether the simple modules in the principal $p$-block of $k\Sym{n}$ are algebraic for $p\geq 7$.

Finally, in Section \ref{sec:sporadic} we consider nine of the sporadic groups: the five Mathieu groups, $HS$, $J_2$, $Suz$ and $He$. In some cases, for example, $HS$ in characteristic $3$, all simple modules are algebraic, and we refer to that section for the specific results that we achieve in this direction.

\medskip

The organization of this paper is as follows: the next section includes the quoted and preliminary results that we need, and Section \ref{sec:SL2q} is concerned with $\SL_2(q)$. The next two sections deal with classical and exceptional groups, and then Lie-type groups in non-defining characteristic in Section \ref{sec:nondefining}. The alternating and symmetric groups are considered in Section \ref{sec:symmetric}, and the final section, Section \ref{sec:sporadic}, is on the sporadic groups.

Throughout this article, all groups are finite and all modules are finite dimensional, unless otherwise specified.

\section{Quoted Results and Preliminaries}

In this section we summarize results gathered from the literature, together with a few other results needed in the sections to come.

The basic properties of algebraic modules may be found in \cite{feit} for example: being algebraic is closed under direct sums, tensor products, taking summands, induction, restriction, Green correspondence and taking sources. This makes the first result trivial.

\begin{prop}\label{prop:p'index} Let $M$ be a $kG$-module, and let $H$ be a $p'$-index subgroup of $G$. We have that $M$ is algebraic if and only if $M\res H$ is algebraic.
\end{prop}

The next three results are in the literature, and we quote them.

\begin{thm}[\cite{craven2009}]\label{thm:abelianSylow2} Let $G$ be a group with abelian Sylow $2$-subgroups, and let $k$ be a field of characteristic $2$. All simple $kG$-modules are algebraic.
\end{thm}

\begin{thm}[{{\cite[Theorem B]{craven2011}}}] \label{thm:nonperiodalg} Let $G$ be a finite group and let $k$ be a field of characteristic $p$. If $M$ is a non-periodic, algebraic $kG$-module, then $\Omega^i(M)$ is non-algebraic for all $i\neq 0$. If $M$ is a periodic, algebraic $kG$-module then $\Omega^i(M)$ is algebraic for all $i\in \Z$.
\end{thm}

\begin{thm}\label{thm:cyclickleinfour} Let $G$ be a finite group, and let $k$ be a field of characteristic $2$. If $B$ is a block of $kG$ with cyclic or Klein-four defect group, then all simple $B$-modules are algebraic.
\end{thm}
\begin{pf} For blocks with cyclic defect group this follows from the simple fact that there are only finitely many indecomposable modules of a group with cyclic vertex, and for Klein-four defect group this result is in \cite{cekl2008}.
\end{pf}

The next easy observation is very useful when dealing with sporadic groups in particular, as it allows us to prove that many modules are non-algebraic almost immediately.

\begin{prop}[$V_4$ restriction test]\label{prop:v4restrictiontest} Let $G$ be a finite group, let $Q$ be a Klein-four subgroup of $G$, and let $k$ be an algebraically closed field of characteristic $2$. If $M$ is a $kG$-module such that $M\res Q$ contains non-trivial odd-dimensional indecomposable summands, then $M$ is non-algebraic.
\end{prop}
\begin{pf} By \cite{conlon1966}, the odd-dimensional indecomposable $kG$-modules are $\Omega^i(k)$ for $i\in\Z$, and $i\neq 0$ for them to be non-trivial. By Theorem \ref{thm:nonperiodalg} these modules are non-algebraic, and hence $M$ is non-algebraic.
\end{pf}

As a remark, by \cite[Theorem 3.4]{archer2008} this result holds with $Q$ replaced by any dihedral $2$-group.

\begin{prop}[\cite{craven2011}]\label{prop:6dimperalg} Let $G=C_p\times C_p$ and let $k=\F_3$. If $M$ is an absolutely indecomposable $kG$-module of dimension either $3$ or $6$, then $M$ is periodic if and only if it is algebraic.
\end{prop}

We remind the reader of \cite[Conjecture E]{craven2011}, which motivates some of the expected results in this article.

\begin{conj}\label{conj:periffalg} Let $k$ be a field of characteristic $p$, and let $G=C_p\times C_p$, where $p$ is an odd prime. If $M$ is an absolutely indecomposable $kG$-module of dimension a multiple of $p$, then $M$ is algebraic if and only if it is periodic.
\end{conj}

Let $T(M)$ denote the (infinite-dimensional) module
\[ T(M)=\bigoplus_{n\geq 1} M^{\otimes n},\]
and let $\mc T(M)$ denote the set of isomorphism classes of indecomposable summands of $T(M)$.

\begin{lem}\label{lem:Mplusalgebraic} Let $G$ be a finite group, and let $k$ be a field of characteristic $p$. Let $M$ be an indecomposable $kG$-module. If $M^{\otimes 2}=M\oplus X$, where $X$ is algebraic, then $M$ is algebraic.
\end{lem}
\begin{pf} Let $\mc X=\mc T(X)$, and let $\mc M$ denote the (finite) set of isomorphism classes of indecomposable summands of $M\otimes A$, as $A$ runs through $\mc X$, together with $M$ itself. We claim that $\mc M=\mc T(M)$. Certainly, as $X$ is a summand of $M^{\otimes 2}$, all modules in $\mc M$ appear in $\mc T(M)$. To see the converse, we claim that if $A\in \mc M$ then all summands of $A\otimes M$ lie in $\mc M$. Since $M\in \mc M$, this will complete the proof. If $A=M$ then $A\otimes M=M\oplus X$, and the claim is true; if $A\in \mc X$ then again clearly all summands of $A\otimes M$ lie in $\mc M$, by construction. If $A$ is a summand of $B\otimes M$, for $B\in \mc X$, then
\[ A\otimes M\mid B\otimes M\otimes M=B\otimes (M\oplus X)=B\otimes M\oplus B\otimes X.\]
The summands of $B\otimes M$ lie in $\mc M$, and the summands of $B\otimes X$ lie in $\mc X\subs \mc M$, so the proof is complete.
\end{pf}

\begin{prop}\label{prop:heartprojcover} Let $p$ be an odd prime, let $k$ be a field of characteristic $p$, and let $G$ be a finite group with non-cyclic Sylow $p$-subgroups. The heart of $\proj(k)$, $\rad(\proj(k))/\soc(\proj(k))$, is non-algebraic.
\end{prop}
\begin{pf} Since $G$ has non-cyclic Sylow $p$-subgroups, $k$ is non-periodic. Since $p$ is odd, $k$ lies on a component $\Gamma$ of the Auslander--Reiten quiver of type $A_\infty$. Write $M_1=k$, and $M_{2i+1}$ for $3\leq 2i+1\leq p$ for the indecomposable modules lying directly above $M_1$ on $\Gamma$ (with $M_3$ the closest to $M_1$). The modules $\Omega^i(k)$ have dimension $\pm 1\bmod p$, and it is easy to see that the modules $\Omega^i(M_j)$ have dimension $\pm j\bmod p$; in particular, $M_p$ has dimension a multiple of $p$. Write $\Gamma'$ for the component of the Auslander--Reiten quiver containing $\Omega(k)$: write $M_2$ for the middle term of the almost-split sequence starting in $\Omega(k)$, and $M_{2i}$ for $4\leq 2i\leq p-1$ for the modules lying directly above $M_2$ on $\Gamma'$.

Let $\mc A_i$ denote the almost-split sequence starting in $\Omega(M_i)$, which has middle term $M_{i-1}\oplus M_{i+1}$. By \cite[Theorem 3.6]{auslandercarlson1986}, taking the tensor product of $\mc A_1$ by $M_i$ yields (up to projectives) $\mc A_i$, and so taking middle terms we get
\[ M_{i-1}\oplus M_{i+1}=\Omega^0(M_i\otimes M_2),\]
for $1<i<p$. Finally, since $M_p$ has dimension a multiple of $p$, $\mc A_1\otimes M_p$ is split by \cite[Theorem 3.6]{auslandercarlson1986} again, and so (modulo projectives) $M_p\otimes M_2=\Omega(M_p)\oplus \Omega^{-1}(M_p)$.

The tensor powers of $M_2$ contain all $M_i$ for $2\leq i\leq p$, and hence contains the summands of $M_2\otimes M_p$. As not both of $M_p$ and $\Omega(M_p)$ can be algebraic, Theorem \ref{thm:nonperiodalg} implies that $M_2$ is non-algebraic, as claimed.
\end{pf}

\begin{prop}\label{prop:MEnonalgebraic} Let $p$ be an odd prime, let $k$ be a field of characteristic $p$, and let $G$ be a finite group with non-cyclic Sylow $p$-subgroups. Write $E$ for the heart of $\proj(k)$.
\begin{enumerate}
\item If $M$ is a non-algebraic module then $E\otimes M$ is non-algebraic.
\item If $M$ is an algebraic module of $p'$-dimension then $E\otimes M$ contains a non-algebraic summand of dimension prime to $p$.
\end{enumerate}
In particular, if $p\nmid\dim M$ then $M\otimes E$ is non-algebraic.
\end{prop}
\begin{pf} If $M$ is non-algebraic then, since $E$ is self-dual and of dimension prime to $p$, $M^{\otimes 2}$, which is non-algebraic, is a summand of $(M\otimes E)^{\otimes 2}$. Hence $(M\otimes E)^{\otimes 2}$, so $M\otimes E$, is non-algebraic, proving (i).

For (ii), notice that, since $M$ is algebraic, if all $p'$-dimensional summands of $E\otimes M$ were algebraic, then the same would be true for $E\otimes M\otimes M^*$ as, if $X$ and $Y$ are absolutely indecomposable modules with one of $\dim X$ and $\dim Y$ a multiple of $p$, the same is true for all summands of $X\otimes Y$ by \cite[Proposition 2.2]{bensoncarlson1986}. However, as $k\mid M\otimes M^*$, we see that $E\mid E\otimes M\otimes M^*$, a contradiction as $E$ is non-algebraic. This proves (ii).
\end{pf}

Using this proposition we may easily prove Proposition \ref{prop:secondrow} from the introduction. Let $M$ be a $kG$-module that lies on the second row of its Auslander--Reiten quiver, so that the almost-split sequence with middle term $M$ is
\[ 0\to \Omega^{-1}(N)\to M\oplus X\to \Omega(N)\to 0\]
where $X$ is zero or a projective module. If $p\nmid \dim M$ then, by \cite[Theorem 3.6]{auslandercarlson1986}, if $E$ denotes the heart of $\proj(k)$, we have that $N\otimes E=M\oplus X$; by Proposition \ref{prop:MEnonalgebraic}, this module, hence $M$, is non-algebraic, as claimed.

\section{$\SL_2(q)$}
\label{sec:SL2q}
In this section, $p$ is an odd prime, $q=p^n$, $k=\F_q$, and $G=\SL_2(q)$. Let $\b G=\SL_2(\bar k)$. We embed $G$ inside $\b G$ in the standard way, as the fixed points under the map $F^n$, where $F:\b G\to\b G$ is the Frobenius map raising each matrix entry to the $p$th power.

It is well known that the simple $k\b G$-modules $L(\lambda)$ are labelled by non-negative integers $\lambda$, and the simple $kG$-modules are the restrictions of $L(\lambda)$ for $0\leq\lambda\leq q$. Let $V_i=L(\lambda)|_G$ for $0\leq i\leq p-1$, the so-called \emph{fundamental modules}. By Steinberg's tensor product theorem, for $\lambda\geq 0$, writing $\lambda=\sum_{i=0}^d j_i p^i$ (for some $d$), we have that $L(\lambda)$ is the tensor product
\[ L(\lambda)\cong\bigotimes_{i=0}^{d} L(\lambda_i)^{\sigma^i}\]
for $0\leq \lambda_i\leq p-1$, where $\sigma$ is the map on the module category obtained by applying the Frobenius morphism (so that it sends $L(\lambda)$ to $L(p\lambda)$, for example). This tensor product theorem obviously restricts to $G$, with $d$ replaced by $n-1$.

We also require tilting modules for $\SL_2(\bar k)$. We will not repeat the definition here, but instead give the properties of them that we need (see \cite{donkin1993}). The indecomposable tilting modules, $T(\lambda)$, are parameterized by non-negative integers, with $T(\lambda)$ having composition factors $L(\lambda)$ and $L(\mu)$ for $\mu<\lambda$, and a tilting module is a sum of the $T(\lambda)$ for various $\lambda$. The tensor product of two tilting modules is a tilting module.

We determine some specific tilting modules now, using \cite[Lemma 5]{erdmannhenke2002}. For $0\leq\lambda\leq p-1$ we have $T(\lambda)=L(\lambda)$; for $p\leq \lambda\leq 2p-2$, writing $\lambda=p+\mu$, $T(\lambda)$ is uniserial of length $3$, with radical layers $L(p-2-\mu)$, $L(\lambda)=L(\mu)\otimes L(1)^\sigma$ and $L(p-2-\mu)$; the module $T(2p-1)$ is simple, and is  $L(2p-1)=L(p-1)\otimes L(1)^\sigma$. Finally, notice that, since $T(\lambda)$ has composition factors $L(\lambda)$ and $L(\mu)$ for $\mu<\lambda$, a tilting module is determined up to isomorphism by its composition factors.

We now identify the modules $V_i$: $V_0$ is the trivial module, $V_1$ is the $2$-dimensional natural module, and $V_i=S^i(V_1)$ is the $i$th symmetric power of $V_i$ for $2\leq i\leq p-1$. As $S^i(M)$ is a summand of $M^{\otimes i}$ for $i\leq p$, we see that if $V_1$ is algebraic then so are all fundamental modules. Furthermore, since $M$ is algebraic if and only if the Frobenius twist $M^\sigma$ is algebraic, we see that if $V_1$ is algebraic then all simple $kG$-modules are algebraic, by Steinberg's tensor product theorem.


The tensor products of $L(\lambda)$ and $L(\mu)$ for $0\leq\lambda\leq\mu\leq p-1$ are easy to describe.

\begin{lem}\label{lem:tensoroffundamentals} Let $0\leq\mu\leq\lambda\leq p-1$.
\begin{enumerate}
\item If $\lambda+\mu\leq p-1$, then
\[ L(\lambda)\otimes L(\mu)=L(\lambda-\mu)\oplus L(\lambda-\mu+2)\oplus\cdots \oplus L(\lambda+\mu).\]
\item If $\lambda+\mu> p-1$ and $\lambda<p-1$ then
\begin{align*} L(\lambda)\otimes L(\mu)&=L(\lambda-\mu)\oplus L(\lambda-\mu+2)\oplus \cdots \oplus L(a)
\\ &\qquad\quad\oplus \begin{cases} L(p-1)\oplus T(p+1)\oplus T(p+3)\oplus \cdots \oplus T(\lambda+\mu) & \mu\text{ even}\\T(p)\oplus T(p+2)\oplus \cdots \oplus T(\lambda+\mu) & \mu\text{ odd}\end{cases},\end{align*}
where $a=2p-(\lambda+\mu+4)$.
\item If $\lambda=p-1$ then
\[ L(\mu)\otimes L(p-1)=\begin{cases} L(p-1)\oplus T(p+1)\oplus T(p+3)\oplus \cdots \oplus T(p+\mu-1) & \mu\text{ even}\\T(p)\oplus T(p+2)\oplus \cdots \oplus T(p+\mu-1) & \mu\text{ odd}\end{cases}.\]
\end{enumerate}
\end{lem}
\begin{pf} Notice that the decompositions suggested in the lemma are into tilting modules, so it suffices to prove that the composition factors of the tensor product $L(\lambda)\otimes L(\mu)$ match those of the decomposition. For (i) and (ii), this is \cite[Lemma 4]{erdmannhenke2002}.

Now suppose that $\mu=1$ and $\lambda=p-1$. By \cite[Theorem 3.2]{dotyhenke2005}, $L(1)\otimes L(p-1)=T(p)$, agreeing with (iii). From this one can easily determine the composition factors of $L(\mu)\otimes L(p-1)$, and this gives (iii).
\end{pf}

\begin{lem}\label{lem:threeistwo} We have
\begin{align*} T(p)\otimes L(p-1)&\cong 2\cdot(T(p)\oplus T(p+2)\oplus \cdots \oplus T(2p-3))\oplus L(2p-1)
\\ &\cong 2\cdot L(p-1)\otimes L(p-2)\oplus L(p-1)\otimes L(1)^\sigma.\end{align*}
Consequently, for $0\leq \lambda,\mu\leq p-1$, we may express $L(1)\otimes L(\lambda)\otimes L(\mu)$ as a sum $A\oplus B\otimes L(1)^\sigma$, where $A$ and $B$ are sums of tensor products $L(i)\otimes L(j)$ for $i,j\leq p-1$, and $B=0$ unless $\lambda=\mu=p-1$, in which case $B=L(p-1)$.
\end{lem}
\begin{pf} Firstly, notice that by Lemma \ref{lem:tensoroffundamentals} and the fact that $L(2p-1)=L(p-1)\otimes L(1)^\sigma$, the second and third expressions are equal. Next, since $L(2p-1)=T(2p-1)$, all expressions are of tilting modules, so it suffices to check that the composition factors coincide. However, the composition factors of $T(p)$ are $L(1)^\sigma$ and $L(p-2)$ twice, so the composition factors of $T(p)\otimes L(p-1)$ are those of $L(p-1)\otimes L(1)^\sigma$ and $L(p-1)\otimes L(p-2)$ twice, which are clearly the same composition factors as the third expression in the lemma.

To see the consequence, note that if $0<\lambda<p-1$ then $L(1)\otimes L(\lambda)\otimes L(\mu)=L(\lambda-1)\otimes L(\mu)\oplus L(\lambda+1)\otimes L(\mu)$, and similarly for $\mu$, so that $\lambda=\mu=p-1$, and we have the displayed equation in the lemma.
\end{pf}

\begin{thm} The module $V_1$ is algebraic.
\end{thm}
\begin{pf} In the case where $d=1$, the Sylow $p$-subgroups of $G$ are cyclic, so the result holds: hence we assume that $d>1$.

We claim that every module in $\mc T(V_1)$ appears as an indecomposable summand of $M\otimes N$ for simple $kG$-modules $M$ and $N$; let $\mc M$ denote the set of such summands, and notice that $X\in\mc M$ if and only if $X^\sigma\in\mc M$. It suffices to show that if $M$ and $N$ are simple $kG$-modules then every summand of $V_1\otimes M\otimes N$ appears in $\mc M$. We proceed by induction on $\dim(M\otimes N)$, noting that if this dimension is $1$ then the result is trivial.

By Steinberg's tensor product theorem, we may write $M\otimes N$ as
\[ Y=\bigotimes_{i=0}^{n-1} (V_{i,1}\otimes V_{i,2})^{\sigma^i}.\]
Write $Y=V_{1,1}\otimes V_{1,2}\otimes Y'$, and consider $V_1\otimes Y$: if $V_{1,1}$ and $V_{1,2}$ are not both $V_{p-1}$ then $V_1\otimes V_{1,2}\otimes V_{1,2}$ is a sum of tensor products of two fundamental modules, by Lemma \ref{lem:threeistwo}, and hence $V_1\otimes Y$ is a sum of tensor products of two simple modules, as required. Hence both $V_{1,1}$ and $V_{1,2}$ are $V_{p-1}$, and $V_1\otimes V_{1,1}\otimes V_{1,2}$ can be written as the sum $A\oplus B\otimes V_1^\sigma$, where $A$ and $B$ are non-zero sums of tensor products of fundamental simple modules. Hence $V_1\otimes Y=A\otimes Y'\oplus V_1^\sigma\otimes B\otimes Y'$; the module $A\otimes Y'$ is a sum of tensor products of two simple modules, so all summands of it lie in $\mc M$, and $(V_1^\sigma\otimes B\otimes Y')^{\sigma^{-1}}$ is a tensor product $V_1\otimes M'\otimes N'$ (as $B=V_{p-1})$) with $\dim(M'\otimes N')<\dim(M\otimes N)$, whence all summands of it lie in $\mc M$ by induction. As $X\in \mc M$ if and only if $X^\sigma\in \mc M$, all summands of $V_1^\sigma\otimes B\otimes Y'$ lie in $\mc M$, as needed.
\end{pf}

\section{Natural Modules for the Classical Groups}\label{classicalgroups}

We begin by recalling Jennings's theorem on the group algebras of $p$-groups. Let $P$ be a finite $p$-group. Define the \emph{dimension subgroups}
\[ \Delta_i(P)=[P,\Delta_{i-1}(P)]\Delta_{\ceiling{i/p}}(P)^p.\]
This is the fastest-decreasing central series whose quotients $\Delta_i(P)/\Delta_{i+1}(P)$ are elementary abelian $p$-groups.

\begin{thm}[Jennings, \cite{jennings1941}] Let $P$ be a $p$-group and $k=\F_p$. Denote by $kP$ the group algebra of $P$ over $k$.
\begin{enumerate}
\item Let $A_i(P)$ be defined by
\[ A_i(P)=\{g\in P\,:\,g-1\in \rad^i(kP)\}.\]
Then $A_i(P)=\Delta_i(P)$.
\item Suppose that we choose $x_{i,j}\in P$ such that $x_{i,j}\Delta_{i+1}$ form a basis of $\Delta_i/\Delta_{i+1}$. Write $X_{i,j}=x_{i,j}-1$. Then
\[ \prod_{i,j} X_{i,j}^{\alpha_{i,j}},\qquad 0\leq \alpha_{i,j}\leq p-1\]
generate $kP$. Furthermore, if the weight of such a product is defined to be $\sum_{i,j} i\alpha_{i,j}$, then all products of weight $i$ form a basis of $\rad^{i-1}(kP)/\rad^i(kP)$, and all products of weight at most $i$ form a basis for $kP/\rad^i(kP)$.
\end{enumerate}
\end{thm}

Now let $P=C_p\times C_p$ be the elementary abelian group of order $p^2$, and let $k=\F_p$. Write $M_i=kP/\rad^i(kP)$. Then Jennings' theorem immediately implies the following result.

\begin{prop}\label{firstprops}\begin{enumerate}
\item The module $kP$ has $2p-1$ radical layers.
\item The module $M_i$ has dimension $i(i+1)/2$ if $i\leq p$.
\item The module $M_i$ is spanned by all monomials in $X$ and $Y$ of degree at most $i-1$, for $i\leq p$.
\end{enumerate}
\end{prop}

In particular, (iii) of this proposition implies the next lemma.

\begin{lem}\label{secondprop} Let $1\leq i\leq p-1$ be an integer. Then $S^i(M_2)=M_{i+1}$.
\end{lem}
\begin{pf} This is obvious if one remembers that $S^i(M_2)$ is spanned by all monomials of degree $i$ in the basis elements of $M_2$, which are $1$, $X$ and $Y$. Thus $S^i(M_2)$ is spanned by all monomials in $X$ and $Y$ of degree at most $i$.
\end{pf}

Finally, recall that if $i<p$ and $k$ is a field of characteristic $p$, then for any $kG$-module $M$ the module $S^i(M)$ is a summand of $M^{\otimes i}$.

\begin{prop}\label{prop:cpcpdim3} The 3-dimensional $kP$-module $M_2$ is non-algebraic.
\end{prop}
\begin{pf} Firstly, notice that both $S_1=S^{p-2}(M_2)$ and $S_2=S^{p-1}(M_2)$ are summands of $M_2^{\otimes (p-2)}$ and $M_2^{\otimes (p-1)}$ respectively, and hence if at least one of $S_1$ and $S_2$ is non-algebraic, then $M_2$ is non-algebraic and the proposition follows.

To see that not both of $S_1$ and $S_2$ are algebraic, we simply note that
\[ S_2=\Omega(S_1)^*.\]
Thus by Theorem \ref{thm:nonperiodalg}, either $S_1$ or $S_2$, and consequently $M_2$, is non-algebraic.
\end{pf}

Using this, we easily deal with the special linear groups.

\begin{prop}\label{prop:sln} Let $q$ be a power of a prime $p$, and let $G$ be the group $\SL_n(q)$, where $n\geq 3$. The natural $kG$-module is non-algebraic.
\end{prop}
\begin{pf} The natural module for $\SL_n(q)$ restricts to the subgroup $\SL_n(p)$ as the natural module for this group, so it suffices to prove the result for this group. Since the natural module for $\SL_n(p)$ is algebraic if and only if the natural module for $\GL_n(p)$ is algebraic (by Proposition \ref{prop:p'index}), it suffices to find a non-algebraic module of dimension $n$ over $\F_p$ for some finite group. This is assured by Proposition \ref{prop:cpcpdim3}, by taking the sum of a $3$-dimensional non-algebraic module with $n-3$ copies of the trivial module for the group $C_p\times C_p$.
\end{pf}

Before we deal with the other classical groups in general, we need to prove a result about $\Alt8=\Omega_6^+(2)$.

\begin{lem}\label{lem:o6+2} Let $k$ be a field of characteristic $2$ and let $G=\Alt8=\Omega_6^+(2)$. The $6$-dimensional natural $kG$-module $M$ (viewing $G$ as $\Omega_6^+(2)$) is non-algebraic.
\end{lem}
\begin{pf} Consider the $8$-point natural permutation representation of $\Alt8$: this is easily seen to be uniserial, with radical layers $k$, $M$ and $k$. Let $H$ denote a (transitive) subgroup of $\Alt8$ isomorphic with $\SL_3(2)=\PSL_2(7)$, acting on the eight points of the projective line. The restriction of $M$ to $H$ is the heart of the permutation module of $\PSL_2(7)=\SL_3(2)$ acting on the projective line, and this is the sum of the two $3$-dimensional simple modules. Since these are both non-algebraic by Proposition \ref{prop:sln}, we see that $M$ is non-algebraic, as claimed.
\end{pf}

(As a remark, it can easily be shown that all non-trivial simple modules in the principal $2$-block of $\Alt8$ are non-algebraic.)

Using the well-known theory of alternating forms (see \cite[Chapter 3]{wilsonrob} for example) it is easy to see that for $n\geq 4$ there is a subgroup $H$ of $G=\Sp_{2n}(q)$ isomorphic with $\Sp_6(q)$, such that the natural module for $G$ restricts to the sum of the natural module for $H$ and a $(2n-6)$-dimensional trivial module. For $G=\Sp_6(q)$, one may embed $H=\SL_3(q)$ as block-diagonal matrices, with $h\in H$ as the top $3\times 3$-matrix. In this case, the restriction of the natural module for $G$ is the sum of the natural module for $H$ and its dual, so that the natural module for $\Sp_{2n}(q)$ is non-algebraic.

The case of orthogonal groups is similar: for $G=\Omega^\pm_{2n}(q)$ with $n\geq 4$ and $q$ even, there is a copy of $H=\Omega_6^+(2)$ embedded so that the natural module for $G$ restricts to the sum of the natural module for $H$ and a $(2n-6)$-dimensional trivial module. For odd $q$, we reduce to $\Omega^+_6(q)$ in the same way as for even $q$, and here again see a diagonally embedded copy of $\SL_3(q)$ in the same way. Hence for orthogonal groups $\Omega^\pm_{2n}(q)$, $n\geq 4$, the natural module is non-algebraic. The case of $\Omega_{2n+1}(q)$, $n\geq 3$, is very similar, and is omitted.

For unitary groups $\SU_n(q)$, things are similar, but only when $n\geq 6$. Consider $\SU_n(q)$ as those matrices $A$ in $\SL_n(q^2)$ such that $(\bar A)^t=A^{-1}$, where $\bar A$ denotes the matrix obtained from $A$ by raising each entry to the $q$th power. One notes that $\SU_6(q)\cap\SL_6(q)$ is the set of all matrices satisfying $A^t=A^{-1}$, so that it is orthogonal for odd $q$ and symplectic for even $q$. As we showed that, for both of these, the natural module is non-algebraic, this completes the statement.

Collating these statements, we have the following.

\begin{prop} Suppose that $k$ is an algebraically closed field of characteristic $p$, let $q=p^a$, and let $G=\Sp_{2n}(q)$, $\SU_{n+3}(q)$, $\SO_{2n}^+(q)$ or $\SO_{2n+2}^-(q)$ ($q$ odd or even), or $\SO_{2n+1}(q)$ ($q$ odd), where $n\geq 3$. The natural $kG$-module is non-algebraic.
\end{prop}

In the case of $\Sp_4(q)$ it is not known in general what happens. For $q=2$, $G=\Sp_4(q)\cong \Sym{6}$, and in this case it is true that all simple $kG$-modules are algebraic, as we shall see in Theorem \ref{thm:altsymmetricweight2}. However, in general it is not clear whether the natural module is algebraic. (Note that the Suzuki groups $\Suz(2^{2n+1})$ are likely to follow the same pattern as $\Sp_4(2^n)$.)

For unitary groups, we are interested in $\SU_n(q)$ for $n=3,4,5$. The group $\SU_3(2)$ is soluble, so the natural module is algebraic, and the natural module for $\SU_3(3)$ is algebraic, by a computer-based proof. For $q\geq 4$ however, it is not clear whether the natural module for $\SU_3(q)$ is algebraic, even using a computer.

For $n=4$, the natural module for $\SU_4(2)$ is non-algebraic, as there is a conjugacy class of $V_4$ subgroups such that the natural module restricts to the sum of a $1$- and $3$-dimensional module, and so is non-algebraic by the $V_4$ restriction test (Proposition \ref{prop:v4restrictiontest}). By embedding $\SU_4(2)$ into $\SU_n(q)$, where $q$ is even and $n=4,5$, we see that the natural module is non-algebraic for even $q$ and all $n\geq 4$. It appears as though this statement also holds for odd $q$, but as of yet the author does not have a proof.

\section{The Natural Module for the Exceptional Groups}
\label{sec:exceptional}

In this section we will prove that the natural module for the exceptional groups of Lie type is non-algebraic in all cases apart from the Suzuki groups. The strategy in the larger-rank cases is as follows: for an algebraic group $\b G$, find a subgroup $\b H$ such that the natural module, restricted to $\b H$, has a simple summand that is known to be non-algebraic by earlier results. For small-rank groups, we use some direct computation and some general results on restricting modules. We use the notation from \cite{liebeckseitz1996} for the weights of algebraic groups.

\begin{prop} If $G$ is the group $^2G_2(q)$, where $q$ is an odd power of $3$, then the $7$-dimensional natural module is non-algebraic.\end{prop}
\begin{pf} Let $G$ be the group $^2G_2(3)=\SL_2(8)\rtimes C_3$, and let $M$ denote the natural $kG$-module, where $k$ is a splitting field of characteristic $3$. This group has a $9$-point permutation representation, and the corresponding permutation module $X$ is uniserial, with radical layers $k$, $M$ and $k$.

Let $P$ be a Sylow $3$-subgroup of $G$. It is generated by an element $x$ of order $9$, which generates a Sylow $3$-subgroup of $\SL_2(8)$, and $y$ of order $3$, acting non-trivially on $\gen x$. Let $Q=\gen{x^3,y}$, which must act transitively on the nine points; thus we have $X\res Q=kQ$, and so $M$ restricts to the heart of $kQ$, which is non-algebraic by Proposition \ref{prop:heartprojcover}. Hence $M$ is non-algebraic.

Clearly the natural module for $^2G_2(q)$ restricts to $^G_2(3)$ as $M$, and so we have proved our result.
\end{pf}

Having dealt with the small Ree groups, we turn our attention to the big Ree groups. The smallest is the Tits group, $^2F_4(2)'$, which will require a computer to analyze. Although there is a subgroup of the Tits group isomorphic with $\PSL_3(3)$, and the natural module restricts to the $26$-dimensional simple module for this group, we will prove in Proposition \ref{prop:psl3q} that this module is, in fact, algebraic. Using a computer, we confirm that there are ten conjugacy classes of subgroups of $G={}^2F_4(2)'$ isomorphic with $V_4$, three of which have a normalizer of order $256$. One of these is important, and we will need it in the next proposition.

\begin{prop} If $G$ is the group $^2F_4(q)'$, where $q$ is an odd power of $2$, and $M$ is the $26$-dimensional natural module for $G$, then $M$ is not algebraic.\end{prop}
\begin{pf} There is a conjugacy class of $G={}^2F_4(2)'$ isomorphic with $V_4$, and with normalizer of order $256$, such that the restriction of $M$ to an element $P$ from this class is
\[ M\res P=4\cdot \proj(k)\oplus \Omega^2(k)\oplus \Omega^{-2}(k).\]
By Proposition \ref{prop:v4restrictiontest}, $M\res P$, and hence $M$, is non-algebraic. Since the natural module for $^2F_4(q)$ restricts to the natural module for $^2F_4(2)'$, we get the result.
\end{pf}

In fact, there are five simple modules in characteristic $2$ for the Tits group $^2F_4(2)'$: three in the principal $2$-block and two projective simple modules. The two non-trivial, non-projective, simple modules -- of dimensions $26$ and $246$ -- are both non-algebraic.

The Suzuki groups are the other twisted groups that only exist in certain circumstances. In this case, nothing is known for any of the groups in characteristic $2$; even if we can prove that (say) the natural $4$-dimensional module for $Sz(8)$ is non-algebraic, then since $Sz(8)$ is not inside $Sz(32)$ we get nothing directly from this.

\bigskip

We now turn to $G_2(q)$, and the case where $q$ is even needs a separate treatment, given that the natural module has dimension $6$ in this case.

\begin{prop}\label{prop:G2q} Let $G$ be the group $G_2(q)$, and let $M$ be the $7$-dimensional natural module for $G$, unless $p=2$, in which case $M$ is the $6$-dimensional natural module. Then $M$ is non-algebraic.
\end{prop}
\begin{pf} Firstly, assume that $q$ is odd. Let $H$ be the subgroup $\SL_3(q)\rtimes C_2$, denoted $K_+$ in \cite[Theorem A]{kleidman1988}. In \cite{kleidman1988}, Kleidman remarks that $M$ restricts to $H$ as the sum of the trivial and a simple $6$-dimensional simple module. It is clear that the restriction of this $6$-dimensional module to $H'\cong \SL_3(q)$ is the direct sum of the two (dual) $3$-dimensional simple modules. By Proposition \ref{prop:sln}, the $3$-dimensional simple modules for $\SL_3(q)$ are non-algebraic, and so $M$ is non-algebraic.

Now suppose that $q$ is even. The group $G_2(2)=\PSU_3(3)\rtimes C_2$ has a conjugacy class of Klein-four subgroups with $63$ elements, and restricting the $6$-dimensional natural module to this class is the sum of two $3$-dimensional indecomposable summands; hence by the $V_4$ restriction test, the natural module for $G_2(2)$ is non-algebraic, and the same therefore holds for $G_2(q)$ by restriction to $G_2(2)$.
\end{pf}

By considering the embedding of $G_2(q)$ inside $^3D_4(q)$, we may prove the next result.

\begin{prop} If $G$ is the group $^3D_4(q)$, and $M$ is the $8$-dimensional natural module for $G$, then $M$ is non-algebraic.
\end{prop}
\begin{pf} Let $H$ be a maximal subgroup of $^3D_4(q)$ isomorphic to $G_2(q)$. If $q$ is odd then the module $M$ restricts as the sum of $k$ and the natural module for $H$ (see the construction in \cite[Section 4.6]{wilsonrob}), so is non-algebraic. For even $q$ this restriction is indecomposable, and one proceeds as in Proposition \ref{prop:G2q}: the conjugacy class of Klein-four subgroups with $2457$ members fulfils the requirement of the $V_4$ restriction test, so the natural module is non-algebraic in this case as well.
\end{pf}

To work with the rest of the groups, we will restrict simple modules for an algebraic group to a subgroup isomorphic with $\b A_2$, so we need some information about such modules. If $\b G$ is an algebraic group, $\b H$ is a subgroup of $\b G$ isomorphic to $\b A_2$, and $M$ is a module for $\b G$ such that the restriction of $M$ to $\b H$ has composition factors only $L(00)$, $L(10)$, $L(01)$ and $L(11)$ (of dimensions $1$, $3$, $3$ and $8$ ($7$ if $p=3$) respectively), then the only non-trivial extension can occur between $L(00)$ and $L(11)$ by the linkage principle, and so if either $L(10)$ or $L(01)$ is a composition factor then it must be a summand. As the restrictions of these modules to a fixed point subgroup $H\cong \SL_3(q)$ is non-algebraic, the module $M$ for the fixed point subgroup $G=\b G^F$ is also non-algebraic.

The next group is $F_4(q)$, which has a $26$-dimensional natural module. This module is simple unless $q$ is a power of $3$, in which case this $26$-dimensional module splits as the sum of a $25$-dimensional module and the trivial module.

\begin{prop} Let $G$ be the group $F_4(q)$, and let $M$ be the module $L(0001)$, which is $26$-dimensional if $q$ is not a power of $3$, and $25$-dimensional if $3\mid q$. Then $M$ is not algebraic.
\end{prop}
\begin{pf} Let $\b G=\b F_4$ and let $\b H\cong \b A_2$ be a Levi subgroup corresponding to the two short roots of the root system: the restriction of $L(0001)$ to this subgroup is a sum of one $L(11)$, three $L(10)$, three $L(01)$, and one $L(00)$ if $p=3$. By the previous remarks, this means that each of the $L(10)$ is a summand of $L(0001)\res H$, and so, taking fixed points under the Frobenius map, there is a subgroup $H\cong \SL_3(q)$ such that the natural module for $H$ is a summand of $M\res H$; hence $M$ is non-algebraic, as claimed.
\end{pf}

We will now deal with the groups $E_6(q)$ and $^2E_6(q)$; there is a $27$-dimensional natural module of this group in all characteristics.

\begin{prop} If $G$ is either $E_6(q)$ or $^2E_6(q)$, and $M$ is the $27$-dimensional natural module, then $M$ is non-algebraic.
\end{prop}
\begin{pf} Let $\b G=\b E_6$, and let $\b H$ be a subgroup of $\b E_6$ isomorphic with $\b F_4$. By \cite[Proposition 2.5]{liebeckseitz1996}, the restriction of the $27$-dimensional natural module for $\b G$ to $\b H$ is the sum of the natural module $L(0001)$ and the trivial module, for all $p\neq 3$. In the case where $p=3$, we use a different subgroup, namely the one isomorphic to $\b C_4$, which only exists for odd primes. By \cite[Proposition 2.5]{liebeckseitz1996} again, the restriction to this $\b C_4$ subgroup is $L(0100)$.

Hence we consider the module $L(0100)$ for $\b G_1=\b C_4$. Let $\b H_1\cong \b A_2$ be the Levi subgroup corresponding to the $A_2$ on the end of the Dynkin diagram. The restriction of $L(0100)$ to this subgroup is a sum of one $L(11)$, three $L(10)$, three $L(01)$, and two $L(00)$ (only one if $p>3$). Hence $L(10)$ is a summand of $L(0100)\res{\b H_1}$.

Taking fixed points under the appropriate Frobenius map, and noting that $F_4(q)$ and $\Sp_8(q)$ ($q$ odd) are contained in $^2E_6(q)$, we see that for all $q$ the module $M$ is non-algebraic since the restriction of $M$ to $F_4(q)$ or $\Sp_8(q)$ is non-algebraic (the latter case because its restriction to the $\SL_3(q)$-Levi is non-algebraic).
\end{pf}

To deal with the groups $E_7(q)$ and $E_8(q)$, we use the fact discussed above for modules for $\b A_2$.

\begin{prop} Let $G$ be either $E_7(q)$ and $E_8(q)$, and let $M$ denote the natural module, of dimension $56$ and $248$ respectively. Then $M$ is non-algebraic.
\end{prop}
\begin{pf} We use the various tables from \cite{liebeckseitz1996}. If a module $N$ for $\b A_2$ has composition factors $L(00)$, $L(10)$, $L(01)$ and $L(11)$, then we say that $N$ has `the desired composition factors'. This will complete the proof for $E_7(q)$ and $E_8(q)$ by taking fixed points under the appropriate Frobenius map.


If $\b G=\b E_7$, then by \cite[Proposition 2.3]{liebeckseitz1996} there is a subgroup isomorphic with a central product $\b A_2\b A_5$, and restricting the $56$-dimensional simple module to the $\b A_2$ factor has the desired composition factors.

Now let $\b G=\b E_8$. By \cite[Proposition 2.1]{liebeckseitz1996}, there is a subgroup that is a central product $\b A_2\b E_6$, and restricting the Lie algebra module to the $\b A_2$ factor has the desired composition factors.
\end{pf}

One should note that the same argument, using \cite[Proposition 2.1]{liebeckseitz1996}, proves that the $kG$-module corresponding to the Lie algebra module is non-algebraic for $\b G$ of type $\b F_4$, $\b E_6$ and $\b E_7$ as well.

\section{Non-Defining Characteristic Groups}
\label{sec:nondefining}
In this short section we use a mixture of theoretical and computational techniques to prove results in characteristics $2$ and $3$. We use known theorems about the fact that the sources of simple modules in non-defining characteristic are `generic', in the sense that they do not depend on the underlying field but only on congruences. We start with the following easy proposition.

\begin{prop}\label{prop:psl2q} Let $G$ be the group $\PSL_2(q)$ for $q$ odd, and let $k$ be a field of characteristic $2$. If $q\equiv -1\bmod 8$ then the non-trivial simple modules in the principal block of $kG$ are non-algebraic, and in all other cases all simple modules are algebraic.
\end{prop}
\begin{pf} We may assume that $\F_4\subs k$, so that $k$ is a splitting field for $G$. If $q\equiv 3,5\bmod 8$, then the Sylow $2$-subgroup of $G$ is Klein-four: hence all simple $kG$-modules are algebraic by Theorem \ref{thm:cyclickleinfour}. If $q\equiv 1\bmod 4$ then by \cite[Theorem 1]{erdmann1977} the source of a non-trivial simple module in $B_0(kG)$ is $2$-dimensional, hence algebraic by Corollary \ref{cor:2dim}. If $q\equiv 3\bmod 4$ then by \cite[Theorem 3]{erdmann1977} the source of a non-trivial simple module in $B_0(kG)$ has dimension $2^{a-1}-1$, where $2^a$ is the order of a Sylow $2$-subgroup of $G$. Since $2^{a-1}-1$ is odd and greater than $1$, and dihedral $2$-groups have no non-trivial odd-dimensional indecomposable modules that are algebraic by \cite[Theorem 3.4]{archer2008}, this module cannot be algebraic, as claimed.

It is easy to see that for $\PSL_2(q)$ there is a unique block with non-cyclic defect group, and so the proof is complete.
\end{pf}

The remaining simple group with dihedral Sylow $2$-subgroups is $\Alt7$. By \cite[Theorem 5]{erdmann1977}, the three simple modules in the principal $2$-block have source either trivial or $2$-dimensional, so these are algebraic, just as in the previous proposition. As the non-principal $2$-block has Klein-four defect group, we see that all simple $k\Alt7$-modules are algebraic, where $k$ is a field of characteristic $2$.

The two odd central extensions of simple groups with dihedral Sylow $2$-subgroups are $3\cdot \Alt6$ and $3\cdot\Alt 7$: in the former case, there is a faithful $2$-block with full defect, and it is splendidly Morita equivalent to the principal $2$-block of $\PSL_2(7)$, so there are non-algebraic simple modules; in the latter case, the faithful $2$-block with full defect has $6$-dimensional and $15$-dimensional simple modules, both of which are algebraic.

Having dealt with simple groups with dihedral Sylow $2$-subgroup, we consider simple groups with semidihedral Sylow $2$-subgroup, in characteristic $2$. Together with $M_{11}$, these are $\PSL_3(q)$ and $\PSU_3(q)$, for certain congruences modulo $4$. In the case where the Sylow $2$-subgroup is semidihedral, work of Erdmann \cite{erdmann1977} has identified the vertices and sources of the simple modules lying in blocks with full defect, allowing us to prove the following proposition.

\begin{prop}\label{prop:psl3q} Let $G$ be a simple group with semidihedral defect group, and let $k$ be a field of characteristic $2$.
\begin{enumerate}
\item If $G\cong \PSL_3(q)$ for $q\equiv 3\bmod 4$, then all simple $kG$-modules in blocks of full defect are algebraic if and only if $q\equiv 3\bmod 8$, and if $q\equiv 7\bmod 8$ then there is exactly one non-algebraic simple module in each block of full defect.
\item If $G\cong \PSU_3(q)$ for $q\equiv 1\bmod 4$, then all simple modules in blocks of full defect are algebraic.
\item If $G\cong M_{11}$ then all simple $kG$-modules are algebraic.
\end{enumerate}
\end{prop}
\begin{pf}
\begin{enumerate}
\item Let $G$ act on the projective plane in the standard way, with $q^2+q+1$ points, and let $H$ denote a point stabilizer. This is of the form $U\rtimes C$, where $C$ is the centralizer of an involution and $U$ is a group of odd order. The group $C$ is a quotient of $\GL_2(q)$ by a central subgroup of odd order: by Proposition \ref{prop:psl2q}, if $q\equiv 7\bmod 8$ then there are non-algebraic simple modules in the principal $2$-block of $\PSL_2(q)$, and the diagonal automorphism amalgamates these two modules, so that the principal $2$-block of $kC$ has a unique non-trivial simple module, and this is non-algebraic.

In \cite[(3,4)]{erdmann1979}, it is shown that the sources of simple modules in blocks of full defect are either trivial source, or uniserial, of dimension $2^{n-3}-1$, where the power of $2$ dividing $|G|$ is $2^n\geq 16$. If $q\equiv 3\bmod 8$ then $n=4$, so all simple modules have trivial source, thus are algebraic. If $q\equiv 7\bmod 8$ then two of the simple modules in the principal block have trivial source, and the other, $M$, has the property that the restriction to $C$ has as a summand the non-trivial simple module in the principal $2$-block of $kC$: this module is non-algebraic, so $M$ is non-algebraic. In \cite[(3,4)]{erdmann1979} it is also stated that each $2$-block of $kG$ contains a simple module whose source is the same as $M$, completing the proof.

\item This follows immediately from \cite[(4.10)]{erdmann1979}, where it is proved that all sources of simple modules in blocks of full defect are of dimension at most $2$, so are algebraic by Corollary \ref{cor:2dim}.

\item This is proved in Proposition \ref{prop:m11}.
\end{enumerate}
\end{pf}

Let $G$ be the group $\SL_3(q)$ for $q\equiv 3\bmod4$. The centralizer of an involution is the group $\GL_2(q)$, and the defect groups of blocks of $\GL_2(q)$ are easy to understand: they are Sylow $2$-subgroups of direct products of $\GL_d(q)$ for $d\leq 2$, so are either of full defect or abelian. Hence all $2$-blocks of $G$ are either of full defect, or of Klein-four or cyclic defect group. Since all simple modules in blocks with the latter two defect groups are algebraic, this means that all simple modules for $\PSL_3(q)$ belonging to $2$-blocks that are not of full defect are algebraic; in particular, all simple modules for $\PSL_3(q)$ are algebraic if $q\equiv 3\bmod 8$.

Similarly, if $G$ is the group $\SU_3(q)$ for $q\equiv 1\bmod 4$ then the centralizer of an involution is $\GU_2(q)$, and again the defect groups are either a Sylow $2$-subgroup of $G$, Klein-four, or cyclic. Hence again, all simple modules lying in $2$-blocks of $\PSU_3(q)$ (for $q\equiv 1\bmod 4$) that are not of full defect are algebraic.

\medskip

We now move on to characteristic $3$, and consider simple groups with Sylow $3$-subgroup $C_3\times C_3$.

\begin{prop}\label{prop:indivgroups3} Let $k$ be a field of characteristic $3$. Let $G$ be one of the groups $\PSL_3(4)$, $\PSU_3(5)$, $\PSL_4(2)$, $\PSL_5(2)$, $\PSU_4(4)$, $\PSU_5(4)$, $\PSp_4(2)$, $\PSp_4(4)$, $\PSL_2(q)$ for $q$ a power of $3$, $\Alt7$, $M_{22}$ and $HS$. All simple modules in $B_0(kG)$ are algebraic.
\end{prop}
\begin{pf} If $G$ is one of $\PSp_4(4)$, $\PSU_4(4)$ and $\PSU_5(4)$, then by \cite{puig1990} there is a splendid Morita equivalence between $B_0(kG)$ and $(C_3\times C_3)\rtimes D_8$, so that all simple $B_0(kG)$-modules are algebraic. Similarly, if $G=\PSL_5(2)$ then there is a splendid Morita equivalence between $B_0(kG)$ and $(C_3\times C_3)\rtimes D_8$ by \cite{koshitanimiyachi2000}, so that all simple $B_0(kG)$-modules are algebraic.

If $G$ is one of $\PSp_4(2)=\Sym6$, $\Alt7$ or $\PSL_4(2)=\Alt8$ then all simple $B_0(kG)$-modules are algebraic by Theorem \ref{thm:altsymmetricweight2} below. Finally, if $G$ is $M_{22}$ or $HS$, then all simple $B_0(kG)$-modules are algebraic by Propositions \ref{prop:m22} and \ref{prop:hs} respectively.
\end{pf}

Notice that the two sporadic groups $M_{11}$ and $M_{23}$ are not on this list: by Propositions \ref{prop:m11} and \ref{prop:m23} they have non-algebraic simple modules.

We are now in a position to state the main theorem about these groups.

\begin{thm} Let $G$ be a finite group, and suppose that a Sylow $3$-subgroup $P$ of $G$ has order at most $9$. Let $k$ be a field of characteristic $3$. All simple $B_0(kG)$-modules are algebraic if and only if $G$ does not have a composition factor $M_{11}$ or $M_{23}$.
\end{thm}
\begin{pf} Let $G$ be a finite group of this form. By \cite[Section 5]{fongharris1993}, there are normal subgroups $H$ and $L$ of $G$ with $L\leq H$, such that $G/H$ and $L$ are $3'$-groups and $H/L$ is a direct product of simple groups and an abelian $3$-group. As $L$ is a normal $3'$-group, the simple $kG$-modules and the simple $k(G/L)$-modules in the principal $3$-block are the same, so we may assume that $L=1$. Also, since $G/H$ is a $3'$-group, a simple $kH$-module $M$ is algebraic if and only if $M\ind G$ is algebraic. Since the simple modules in $B_0(kG)$ are summands of $M\ind G$ for $M$ a simple $B_0(kH)$-module, we may assume that $G=H$, in which case $G$ is a direct product of simple groups with abelian Sylow $3$-subgroups and an abelian $3$-group. As $|P|=9$, if $G$ is the direct product of more than one group then the factors have cyclic Sylow $3$-subgroup, whence all simple $kG$-modules are algebraic, and in the other case $G$ is simple. There is a splendid Morita equivalence between $B_0(kG)$ and $B_0(kH)$, where $H$ is one of $\PSL_3(4)$, $\PSU_3(5)$, $\PSL_4(2)$, $\PSL_5(2)$, $\PSU_4(4)$, $\PSU_5(4)$, $\PSp_4(2)$, $\PSp_4(4)$, $\PSL_2(9)$, $\Alt7$, $M_{11}$, $M_{22}$, $M_{23}$, and $HS$ (see \cite{koshitanikunugi2002} and the references contained therein). By hypothesis, $G\not\cong M_{11}$ and $G\not\cong M_{23}$, so by Proposition \ref{prop:indivgroups3} the theorem is proved.
\end{pf}

\section{Alternating and Symmetric Groups}
\label{sec:symmetric}
In this section we consider blocks of alternating and symmetric groups having defect groups $C_p\times C_p$. If $p=2$ then this is the Klein-four group, and all simple modules in such blocks are algebraic by Theorem \ref{thm:cyclickleinfour}. Hence we consider odd primes $p$: in general, knowing the sources of simple modules in these blocks is a difficult problem, but for small primes this can be done.

The primary tool in this section will be the following theorem of the author \cite{craven2010un3}.

\begin{thm} Let $B$ be a $p$-block of a symmetric group, and suppose that $B$ has defect group $C_p\times C_p$. If $M$ is any simple $B$-module then the source of $M$ is isomorphic to the source of a simple module from either $B_0(k\Sym{2p})$ or $B_0(k[\Sym p\wr C_2])$.
\end{thm}

The proof of this theorem uses various combinatorial techniques, and is entirely theoretical. Since all simple $k\Sym p$-modules are algebraic, we see that all simple $k[\Sym p\wr C_2]$-modules are algebraic. Hence we get the following corollary.

\begin{cor} If all simple $k\Sym{2p}$-modules are algebraic, then all simple $B$-modules are algebraic, where $B$ is any $p$-block of any symmetric group with defect group $C_p\times C_p$.
\end{cor}

Thus we need to determine whether the simple $k\Sym{2p}$-modules are algebraic. In general this is unsolved, but for $p=3$ and $p=5$ we can prove that they are.

\begin{thm}\label{thm:altsymmetricweight2} If $p=3$ or $p=5$, and $B$ is a $p$-block of a symmetric or alternating group with defect group $C_p\times C_p$, then all simple $B$-modules are algebraic.
\end{thm}
\begin{pf} It suffices by the previous corollary to prove that result for $\Sym6$ and $p=3$, and $\Sym{10}$ and $p=5$. Let $k$ be a field of characteristic $p$. If $p=3$ then $\Alt6=\PSL_2(9)$, and all simple $k\PSL_2(9)$-modules are algebraic by Theorem \ref{thm:SL2}, so that all simple $k\Sym6$-modules are algebraic, proving the theorem.

Now let $p=5$, and consider $\Alt{10}$. The $8$-dimensional simple $k\Alt{10}$-module has $8$-dimensional source $A_2$, and there are twenty-six different non-projective, indecomposable modules appearing in tensor powers of $A_2$: twenty-two periodic modules, $A_2$, the trivial module $k=A_1$, and the sources $A_3$ and $A_7$ of the simple modules of dimensions $28$ and $56$, which are hence also algebraic.

The $34$-dimensional simple module has a $9$-dimensional source $A_4$, and all nine indecomposable summands of tensor powers of $A_4$ lie in $A_4^{\otimes 2}$. The two $35$-dimensional simple, periodic, modules have $10$-dimensional sources $A_5$ and $A_6$, and there are ten different non-projective indecomposable summands in their tensor powers. (These modules are swapped by the outer automorphism of $\Alt{10}$, so it suffices to check one of them.)

The remaining simple modules in the principal block have dimensions $133$, $133$ and $217$, and the sources have dimensions $8$, $8$ and $42$ respectively. The two $133$-dimensional modules are swapped by the outer automorphism of $\Alt{10}$, and so it suffices to check one of them. There are thirty-four different non-projective summands in $T(A_8)$, with thirty of them periodic.

The last module, $A_{10}$, has dimension $42$. The exterior square $\Lambda^2(A_{10})$ contains a summand $M$ of dimension $47$. It is easier to prove that $M$ is algebraic than $A_{10}$ is, and indeed $M$ is algebraic, with $32$ different indecomposable modules appearing in $T(M)$, of which the largest has dimension $80$. Since $A_{10}^{\otimes 2}$ is a sum of indecomposable modules that also appear in $M^{\otimes 2}$, we see that $A_{10}$ is algebraic, completing the proof for this $5$-block. As $\Alt{10}$ has only one block with non-cyclic defect, all simple $k\Alt{10}$-modules are algebraic, and so all simple $k\Sym{10}$-modules are algebraic, completing the proof.
\end{pf}

\section{Sporadic Groups}
\label{sec:sporadic}
In this section we will examine nine sporadic simple groups: the five Mathieu groups, $HS$, $J_2$, $Suz$ and $He$. Having analyzed the smallest Janko group $J_1$ in \cite{craven2009}, this brings the total to ten.

The results here are almost entirely computer driven, and so we do not provide proofs of our statements that may be checked on the computer. In particular, if a simple module fails the $V_4$ restriction test, we simply state that it does, without going into details, unless it is non-trivial to prove this, even with a computer.

We deal with the Mathieu groups in detail, then only consider certain other sporadic groups, of particular interest, either because they have a block with abelian defect group or because all simple modules are algebraic for some prime $p$ with non-cyclic Sylow $p$-subgroups. More detailed proofs are available in the author's thesis \cite{craventhesis} for these groups in the majority of cases, and we are brief in our justification.

In this section, the simple modules in the principal $p$-block are labelled $k=S_1,S_2,\dots,S_n$, ordered by increasing dimension. In one case we consider a non-principal block explicitly, and its simple modules are denoted $T_1,\dots,T_m$, again ordered by dimension. Information on the simple modules is available in \cite{abc}.

\subsection{The Mathieu Groups}

We deal with the five Mathieu groups here.

\begin{prop}\label{prop:m11} Let $G$ be the simple group $M_{11}$. If $p=2$ then all simple $kG$-modules are algebraic, and if $p=3$ then a simple $kG$-module is algebraic if and only if it is self-dual, so that there are four simple $kG$-modules in the principal block that are not algebraic.
\end{prop}
\begin{pf} Firstly let $p=2$. The subgroup $M_9\rtimes C_2$ has index $55$, and the permutation representation on this subgroup is semisimple, the sum of the three modules in the principal block. Hence all simple modules in $B_0(kG)$ are trivial source, so algebraic. The simple modules outside of the principal block are projective, so are also algebraic, proving the result for $p=2$.

If $p=3$, then the standard $11$-point permutation representation is semisimple, with non-trivial submodule the $10$-dimensional self-dual simple module $S_4$. The sources of the $10$-dimensional simple modules $S_5$ and $S_6$ that are not self-dual are $\Omega^{\pm 2}(k)$, and so are not algebraic. If $M$ is one of the $5$-dimensional simple modules $S_2$ and $S_3$ then $\Lambda^2(M)$ is one of $S_4$ or $S_5$, and as $M^{\otimes 2}=\Lambda^2(M)\oplus S^2(M)$, we see that $M$ is not algebraic. The $24$-dimensional simple module $S_7$ has a $6$-dimensional source, and is algebraic since it is periodic, via Proposition \ref{prop:6dimperalg}. The only simple module not lying in $B_0(kG)$ is projective, and hence algebraic, completing the proof.
\end{pf}

\begin{prop} Let $G$ be the simple group $M_{12}$. If $p=2$ or $p=3$ then the non-trivial simple modules in $B_0(kG)$ are non-algebraic, and the simple modules outside the principal block are algebraic.
\end{prop}
\begin{pf} Firstly suppose that $p=2$. The non-trivial simple modules lying in $B_0(kG)$ fail the $V_4$ restriction test, and the non-principal block has Klein-four defect group, so all simple modules in it are algebraic by Theorem \ref{thm:cyclickleinfour}. Hence the result is proved.

Now let $p=3$. There are two conjugacy classes of subgroups isomorphic with $M_{11}$, with representatives $H_1$ and $H_2$. Restricting the $10$-dimensional simple modules $S_1$ and $S_2$ to the $H_i$ proves that these are non-algebraic, since the restriction of $S_1$ to one of the $H_i$ is the sum of the two $5$-dimensional, non-algebraic, modules for $M_{11}$, and similarly with $S_2$ (and the other of the $H_i$). As $S_1^{\otimes 2}$ is the sum of $k$, the $54$-dimensional projective simple module and the $45$-dimensional simple $S_7$, we see that $S_7$ is non-algebraic (as $S_1^{\otimes 2}$ is non-algebraic, and the other summands are algebraic). Taking $S_2^{\otimes 2}$ yields the other $45$-dimensional simple $S_8$ as a summand of the tensor square, and so it is non-algebraic also.

It remains to discuss the $15$-dimensional simples $S_4$ and $S_5$, and the $34$-dimensional simple $S_6$. We first deal with $S_4$: there are two classes of subgroups of order $9$ with $220$ conjugates, and let $Q$ denote one of these. The restriction of $S_4$ to $Q$ is
\[ S_4\res Q=21\cdot \mc P(k)\oplus M\oplus \Omega^{-1}(M),\]
where $M$ is a $15$-dimensional, non-periodic indecomposable module. (As $S_4$ and $S_5$ are dual, we might need to consider $S_5$ for this decomposition to hold.) As not both of $M$ and $\Omega^{-1}(M)$ can be algebraic, $S_4$ (and hence $S_5$) is not algebraic. The module $S_6$ has vertex a Sylow $3$-subgroup and $7$-dimensional source $S$: the module $S\res Q$ is the heart of the projective indecomposable module for $Q$, and is hence non-algebraic by Proposition \ref{prop:heartprojcover}. This completes the proof.
\end{pf}

\begin{prop}\label{prop:m22} Let $G$ be the simple group $M_{22}$. If $p=2$ then the non-trivial simple modules are non-algebraic, and if $p=3$ then all simple modules are algebraic.
\end{prop}
\begin{pf} Using the $V_4$ restriction test, we get the result for $p=2$. For $p=3$, all non-principal blocks have cyclic defect groups, so it suffices to consider the principal $3$-block. In this case, there is a splendid Morita equivalence with the principal block of the Mathieu group $M_{10}=\Alt6.2_2$, all of whose simple modules are algebraic by Theorem \ref{thm:altsymmetricweight2} and Proposition \ref{prop:p'index}. This completes the proof.
\end{pf}

\begin{prop}\label{prop:m23} Let $G$ be the simple group $M_{23}$. If $p=2$ then all non-trivial simple modules in the principal block, apart from that of dimension $252$, are non-algebraic. If $p=3$ then a simple module is algebraic if and only if it does not have dimension $104$.
\end{prop}
\begin{pf} Using the $V_4$ restriction test, we get the result for $p=2$. For $p=3$, all non-principal blocks have cyclic defect groups, so it suffices to consider the principal $3$-block. All simple modules in $B_0(kG)$, apart from the two $104$-dimensional simple modules $S_3$ and $S_4$, are trivial-source modules, and hence algebraic.

The modules $S_3$ and $S_4=S_3^*$ are non-algebraic: to see this, we note that the source of (say) $S_3$ is a $5$-dimensional module $M$, and $M^{\otimes 2}$ has $\Omega^2(k)$ as a summand. As $\Omega^2(k)$ is non-algebraic, $M$ is non-algebraic, as claimed.
\end{pf}

The $252$-dimensional module has algebraic restriction to all elementary abelian subgroups of $M_{23}$, so it is not clear whether this module is algebraic or not. However, the $252$-dimensional simple module for $M_{24}$ is non-algebraic, and it restricts to this simple module, which suggests that it is probably not algebraic.

\begin{prop}\label{prop:m24} Let $G$ be the simple group $M_{24}$. If $p=2$ then all non-trivial simple modules in the principal block, apart from that of dimension $320$ or $1792$, are non-algebraic. If $p=3$ then a simple module is algebraic if and only if it does not have dimension $770$ (there are two simple modules with this dimension).
\end{prop}
\begin{pf} For $p=2$, the $V_4$ restriction test proves that the non-trivial simple modules, with the exception of the two mentioned, are non-algebraic. Now suppose that $p=3$: all non-principal blocks have cyclic defect group, so it suffices to consider the seven simple modules in the principal block. The $22$-dimensional simple $S_2$ has vertex a Sylow $3$-subgroup $P$ of $G$, and $4$-dimensional source, which has a kernel of order $3$. It is easy to check by direct calculation that this module is algebraic. Since $\Lambda^2(S_2)=S_3$, the simple module of dimension $231$, and $S_2^{\otimes 2}=\Lambda^2(S_2)\oplus S^2(S_2)$, we see that $S_3$ is also algebraic. The $483$-dimensional simple module $S_4$ has trivial source (it is a summand of a permutation module of dimension $759$), and is hence algebraic. Finally for the algebraic modules, the largest simple module $S_7$, of dimension $1243$, has a $19$-dimensional source $M$ (and vertex $P$). The tensor square of this module is (up to projectives) the sum of $k$ and two non-isomorphic $18$-dimensional modules, induced from subgroups of order $9$ in $P$. Hence $S_7$ is algebraic if and only the $6$-dimensional sources $N_1$ and $N_2$ of the $18$-dimensional indecomposable modules are algebraic. The fact that the $N_i$ are periodic, together with Proposition \ref{prop:6dimperalg}, completes the proof that $S_7$ is algebraic.

We turn our attention to the the two (dual) $770$-dimensional simple modules $S_5$ and $S_6$, which are non-algebraic. If $M$ denotes the $5$-dimensional source of $S_5$, then $M^{\otimes 2}$ has as a summand either $\Omega^2(k)$ or $\Omega^{-2}(k)$, so neither $S_5$ nor $S_6=S_5^*$ is algebraic. This completes the proof.
\end{pf}

Again, it is not known whether or not the $320$- and $1792$-dimensional simple modules are algebraic, although the author considers it unlikely.

\subsection{$HS$, $J_2$, $Suz$ and $He$}

Here we consider the Higman--Sims, second Janko, Suzuki and Held sporadic groups. In the case of $HS$, in characteristic $3$ then all simple modules are algebraic, and the same is true of $J_2$ in characteristics $3$ and $5$. For the larger two sporadics however, things are more complicated: it appears possible that all simple modules for $Suz$ are algebraic in characteristic $5$, as at least the smallest six are, and for $He$ \emph{almost} all of the simple modules are algebraic for both $p=3$ and $p=5$.

\begin{prop}\label{prop:hs} Let $G$ be the simple group $HS$. If $p=2$ then a non-trivial simple module is algebraic if and only if it lies outside the principal block. If $p=3$ then all simple modules are algebraic.
\end{prop}
\begin{pf} The result for the prime $2$ follows easily from the $V_4$ restriction test for the principal block and, since the only non-principal block has Klein-four defect group, for those not in the principal block by Theorem \ref{thm:cyclickleinfour}.

Now let $p=3$. There are seven simple modules in each block of $kG$ with non-cyclic defect group, labelled $S_1$ to $S_7$ for the principal block and $T_1$ to $T_7$ for the non-principal block. The modules $S_1$, $S_2$, $S_3$, $T_3$, $T_4$, $T_5$, $T_6$ and $T_7$ are all trivial-source modules, hence algebraic.

There is a subgroup $H\cong M_{22}$ of $G$, of index $100$. All simple $kH$-modules are algebraic by Proposition \ref{prop:m22}, and the restrictions of $S_4$, $T_1$ and $T_2$ to $H$ are semisimple, hence these modules are algebraic by Proposition \ref{prop:p'index}. Also, $S_7\res H$ is the sum of simple and projective modules, so $S_7$ is also algebraic.

This leaves us with $S_5$ and $S_6$. Since $S_6=\Lambda^2(T_1)$, and $T_1^{\otimes 2}=\Lambda^2(T_1)\oplus S^2(T_1)$, this is algebraic. Finally, $S_5$ has a $10$-dimensional source $M$, and $M^{\otimes 2}$ is the sum of a trivial module, nine free modules, and three $6$-dimensional periodic indecomposable modules, which are algebraic by Proposition \ref{prop:6dimperalg}. As $M^{\otimes 2}$ is algebraic, $M$ is algebraic so $S_5$ is algebraic, completing the proof.
\end{pf}

There are unfortunately no results on the (non-trivial) simple modules in the principal $5$-block of $HS$.

\begin{prop} Let $G$ be the Janko group $J_2$. If $p=2$ then all simple modules outside the principal block are algebraic, and all non-trivial simple modules in the principal block are non-algebraic. If $p=3$ or $p=5$ then all simple modules are algebraic.
\end{prop}
\begin{pf} If $p=2$ then the $V_4$ restriction test is enough to prove the result for the principal block and, since the only non-principal block has Klein-four defect group, for those simple modules not in the principal block by Theorem \ref{thm:cyclickleinfour}.

If $p=3$, then all simple modules lying outside the principal block lie in blocks with cyclic defect groups, so that they are all algebraic. There are eight simple modules lying in $B_0(kG)$, labelled as usual $S_1$ to $S_8$. Apart from the trivial module $k=S_1$ and the 133-dimensional module $S_8$, the other six come in pairs, of dimension $13$, $21$ and $57$, each defined over $\F_9$. Because of this, proving (for example) that $S_2$ is algebraic proves that $S_3$ is algebraic. We have firstly that
\[ S_2\otimes S_2=k\oplus S_4\oplus S_6\oplus T_2,\]
where $T_2$ is the $90$-dimensional simple module in the second block. If $S_2$ is algebraic then so are $S_4$ and $S_6$, and by applying the Frobenius morphism on $\F_9$, we see that $S_3$, $S_5$ and $S_7$ are also algebraic. Finally, $S_2\otimes S_3=S_8\oplus T_1$ (where $T_1$ is the $36$-dimensional simple module in the second block), so that $S_8$ is also algebraic. Thus it remains to prove that $S_2$ is algebraic. Let $A_2$ denote the source of $S_2$, which is simply $S_2\res P$. There are nine non-projective modules in $\mc T(A_2)$, so it is algebraic.

Finally, suppose that $p=5$. There are six simple modules in the principal block of $kG$, and all non-principal blocks have cyclic defect groups. The $14$-dimensional simple module $S_2$ has indecomposable restriction to $P$, and hence its source $A_2$ is $14$-dimensional. It is algebraic, and there are twenty-one non-projective modules in $\mc T(A_2)$. Using this we may prove that all simple modules are algebraic: we have
\[ S_2\otimes S_2=k\oplus S_2\oplus S_3\oplus T_1\oplus T_2,\]
where $T_1$ and $T_2$ are the two simple modules lying in the second block (with defect $1$), of dimensions 70 and 90 respectively. Hence $S_3$ is algebraic, and in addition
\[ S_3\otimes S_3=k\oplus S_2\oplus S_3\oplus S_4\oplus S_5\oplus S_6\oplus T_2,\]
so that all simple modules in $B_0(kG)$ are algebraic, as required.
\end{pf}

\begin{prop} Let $G$ be the Suzuki sporadic group $Suz$. If $p=2$ then all non-trivial simple modules in the principal block are non-algebraic. If $p=3$ then all non-trivial simple modules in the principal block are non-algebraic, except possibly for the $8436$- and $32967$-dimensional modules, and all simple modules outside of the principal block are algebraic. If $p=5$ then the smallest six simple modules in the principal block are algebraic.
\end{prop}
\begin{pf} Firstly, let $p=2$. We apply the $V_4$ restriction test, but since the simple modules for $Suz$ are quite large, we describe completely the restrictions involved. Let $Q$ be a representative from the conjugacy class of $V_4$ subgroups with $1216215$ members. We have decompositions of:
\begin{align*}
S_2\res Q&=8\cdot k\oplus \Omega^2(k)\oplus \Omega^{-2}(k)\oplus 20\cdot\proj(k)\oplus X_1;
\\ S_4\res Q&=6\cdot k\oplus 4\cdot \Omega^2(k)\oplus 4\cdot \Omega^{-2}(k)\oplus 24\cdot \proj(k);
\\S_5\res Q&=8\cdot \Omega(k)\oplus 8\cdot \Omega^{-1}(k)\oplus 2\cdot \Omega^3(k)\oplus 2\cdot \Omega^{-2}(k)\oplus 124\cdot \proj(k);
\\S_7\res Q&=4\cdot k\oplus 4\cdot \Omega(k)\oplus 4\cdot \Omega^{-1}(k)\oplus \Omega^4(k)\oplus \Omega^{-4}(k)\oplus 148\cdot \proj(k);
\\S_8\res Q&=8\cdot k\oplus 16\cdot \Omega(k)\oplus 16\cdot \Omega^{-1}(k)\oplus 4\cdot \Omega^2(k)\oplus 4\cdot \Omega^{-2}(k)\oplus 4\cdot \Omega^3(k)\oplus 4\cdot \Omega^{-3}(k)\oplus 808\cdot \proj(k);
\\S_9\res Q&=4\cdot k\oplus 4\cdot \Omega(k)\oplus 4\cdot \Omega^{-1}(k)\oplus \Omega^4(k)\oplus \Omega^{-4}(k)\oplus 1116\cdot \proj(k);
\\S_{10}\res Q&=24\cdot k\oplus 16\cdot \Omega^2(k)\oplus 16\cdot \Omega^{-2}(k)\oplus 4\cdot \Omega^4(k)\oplus 4\cdot \Omega^{-4}(k)\oplus 1184\cdot \proj(k);
\\S_{11}\res Q&=8\cdot k\oplus 6\cdot \Omega(k)\oplus 6\cdot \Omega^{-1}(k)\oplus 6\cdot \Omega^3(k)\oplus 6\cdot \Omega^{-3}(k)\oplus 2264\cdot \proj(k)\oplus X_2;
\\S_{13}\res Q&=12\cdot k\oplus 2\cdot \Omega^2(k)\oplus 2\cdot \Omega^{-2}(k)\oplus 2596\cdot \proj(k)\oplus X_3.\end{align*}
Here $X_1$, $X_2$ and $X_3$ are the sums of six, seventy-two and forty-four $2$-dimensional indecomposable modules respectively. Taking Galois conjugates of the decompositions for $S_2$, $S_{11}$ and $S_{13}$ give the decompositions for $S_3$, $S_{12}$ and $S_{14}$, and taking the dual of the decomposition for $S_5$ gives the decomposition for $S_6$. Therefore all non-trivial simple modules in $B_0(kG)$ are non-algebraic, as claimed.

Now let $p=3$. There are eight conjugacy classes of subgroups of $G$ that are elementary abelian of order $9$, and let $Q$ be a representative from the class with $38438400$ members. Write $E$ for the heart of $\proj(k)$ and $M_1$ for the module $\soc^2(\proj(k))$. There exist $M_2$, a $9$-dimensional, self-dual, non-periodic indecomposable module, and $M_3$, a self-dual, $37$-dimensional indecomposable module, such that the restrictions of the simple modules $S_i$ in the principal block are as follows:
\begin{align*}
S_2\res Q&=E\oplus  \Omega(M_1)\oplus \Omega^{-1}(M_1^*)\oplus 3\cdot \proj(k);
\\ S_3\res Q&=M_2\oplus \Omega^2(M_1)\oplus \Omega^{-2}(M_1^*)\oplus 5\cdot\proj(k);
\\ S_4\res Q&=k\oplus M_2\oplus \Omega^2(M_1)\oplus \Omega^{-1}(M_1^*)\oplus 28\cdot \proj(k);
\\ S_5\res Q&=\Omega(M_1)\oplus \Omega^{-1}(M_1^*)\oplus \Omega^3(M_1^*)\oplus \Omega^{-3}(M_1)\oplus 41\cdot\proj(k);
\\ S_6\res Q&=E\oplus  \Omega(M_1)\oplus \Omega^{-1}(M_1^*)\oplus 68\cdot \proj(k);
\\ S_7\res Q&=M_1\oplus M_1^*\oplus \Omega^3(M_1^*)\oplus \Omega^{-3}(M_1)\oplus \Lambda^2(E)\oplus 209\cdot\proj(k);
\\ S_8\res Q&=\Omega(M_1)\oplus \Omega^{-1}(M_1^*) \oplus \Omega^2(M_1)\oplus \Omega^{-2}(M_1^*)\oplus \Omega(M_2)\oplus \Omega^{-1}(M_2)\oplus \Omega^3(M_1)\oplus \Omega^{-3}(M_1^*)
\\ &\;\;\;\oplus \Omega^4(M_1^*)\oplus \Omega^{-4}(M_1)\oplus 299\cdot \proj(k);
\\ S_9\res Q&=M_2\oplus \Omega^2(M_1)\oplus \Omega^{-1}(M_1^*)\oplus 528\cdot\proj(k);
\\ S_{11}\res Q&=M_1\oplus M_1^*\oplus 2\cdot \Omega(M_1)\oplus 2\cdot \Omega^{-1}(M_1^*)\oplus \Omega^3(M_1^*)\oplus \Omega^{-3}(M_1)\oplus \Omega(M_2)\oplus \Omega^{-1}(M_2)\oplus 1622\cdot \proj(k);
\\ S_{12}\res Q&=2\cdot M_1\oplus 2\cdot M_1^*\oplus \Omega^3(M_1^*)\oplus \Omega^{-3}(M_1) \oplus \Omega^2(k)\oplus \Omega^{-2}(k)\oplus M_3\oplus 2150\cdot\proj(k);
\end{align*}
This proves that all non-trivial simple modules $S_i$ are non-algebraic, apart from $S_{10}$ -- it appears to be non-algebraic, but the proof of this currently eludes the author -- and the $32967$-dimensional simple module $S_{13}$. The latter module has the property that the restriction of it to any conjugacy class of subgroup of order $3$ is free, and so the restriction to any conjugacy class of subgroup of order $9$ is periodic, which makes it difficult to prove that it is non-algebraic; indeed, its restriction to $Q$ is even free. (The size of the module makes a more detailed analysis difficult.)

We turn to the modules outside the principal block. There is only one non-principal block $B$ with non-cyclic defect groups, and in \cite{koshitanikunugiwaki2004}, it is proved that $B$ is splendidly Morita equivalent to the principal $3$-block of $\PSL_3(4)$; hence all simple $B$-modules are algebraic by Proposition \ref{prop:indivgroups3}.

If $p=5$, then the $1001$-dimensional simple module $S_4$ is trivial source, hence algebraic. The $143$-dimensional simple module $S_2$ has a $28$-dimensional source $A_2$, and we have
\[ A_2\otimes A_2=k\oplus A_2\oplus X\oplus 29\cdot \proj(k),\]
where $X$ is the ($30$-dimensional) sum of the permutation modules on the six subgroups of $P$ of order $5$. Hence $A_2$, and $S_2$, are algebraic. The $363$-dimensional simple module $S_3$ has a $13$-dimensional source $A_3$. There are $37$ different non-projective summands of $T(A_3)$, the largest of which has dimension $120$; as this number is finite, $A_3$ is algebraic.

The tensor product of $S_2$ and $S_3$ is a module whose only composition factor in the principal block is $S_6$, and hence this module is algebraic since it must be a summand of $S_2\otimes S_3$. Both the $3289$-dimensional simple module $S_5$ and the $11869$-dimensional simple module $S_6$ have $74$-dimensional sources $A_5$ and $A_6$, and although they are not isomorphic, they are $\Aut(P)$-conjugate. Hence $S_5$ is also algebraic.
\end{pf}


The reason that simple modules of larger dimension for $p=5$ were not considered is simply that their dimensions are too high to be easily constructed on a computer. It is possible that the simple modules of larger dimension are also algebraic. As technology improves more simple modules will be easily dealt with; for example, the $41822$-dimensional simple module is a composition factor of $\Lambda^2(S_3)$, and this module is \emph{almost} able to be easily constructed on a modern computer without more sophisticated techniques, such as the condensation method.

\begin{prop} Let $G$ be the Held sporadic group $He$. If $p=2$ then a simple module is algebraic if and only if it is trivial or lies outside the principal block. If $p=3$ then a simple module is algebraic if and only if it does not have dimension $6172$ or $10879$, and if $p=5$ then the simple modules with dimension $1$, $51$, $104$, $153$, $4116$, $4249$, and $6528$ are algebraic.
\end{prop}
\begin{pf}  First, let $p=2$. The only non-principal block of $kG$ that is not of defect zero has Klein-four defect, and so all simple modules outside the principal block are algebraic. For the modules in the principal block, the $V_4$ restriction test proves the result: for all but the largest two simple modules, of dimension $2449$, a $V_4$ subgroup from one of the two conjugacy classes with $437325$ members will work, whereas for the largest two simple modules, a $V_4$ subgroup from one of the two conjugacy classes with $5247900$ members is necessary.

We next examine $p=3$. The $679$-dimensional simple module $S_2$ has a $4$-dimensional source $A_2$ (which has a kernel of order $3$): there are seven elements in $\mc T(A_2)$, so it is easy to prove that it is algebraic. The $1275$-, $3673$- and $6272$-dimensional simple modules $S_3$, $S_4$ and $S_6$ have trivial source, so are also algebraic.

However, the $6172$-dimensional simple module $S_5$ has a $7$-dimensional source $A_5$ which, when restricted to two of the four maximal subgroups of $P$ (the two that lie in the $G$-conjugacy class of size $2332400$) is the heart of the free module. Hence $A_5$, and so $S_5$, is non-algebraic by Proposition \ref{prop:heartprojcover}. Finally, the $10879$-dimensional has $16$-dimensional source $A_7$, and when restricted to the same maximal subgroups as the previous module, is the sum of a free module and the heart of the free module. Hence $A_7$, and thus $S_7$, is non-algebraic, completing the proof for modules in the principal $3$-block.

In \cite{koshitanikunugiwaki2004}, it is proved that the non-principal $3$-block $B$ of $He$ is splendidly Morita equivalent to its Brauer correspondent, and hence all simple $B$-modules are trivial source, hence algebraic.

Now suppose that $p=5$. The $51$-dimensional simple modules $S_2$ and $S_3$ have trivial source, so are algebraic. The $104$-dimensional simple module $S_4$ has a $29$-dimensional source $A_4$, and there are fifty-eight different indecomposable summands in $T(A_4)$, the largest of which has dimension $129$; thus $S_4$ is algebraic as well. The exterior square of $S_4$ has the $4116$-dimensional simple module $S_{11}$ as a summand, so this module (which has a $66$-dimensional source $A_{11}$) is also algebraic.

The $153$-dimensional simple modules $S_5$ and $S_6$, and the $6528$-dimensional simple module $S_{13}$, all have the same $28$-dimensional source $A_5$, and $A_5^{\otimes 2}$ is the sum of $A_5$, $k$, the six $5$-dimensional permutation modules on subgroups of $P$ of order $5$, and a free module: hence $A_5$ is algebraic by Lemma \ref{lem:Mplusalgebraic}, and so $S_5$, $S_6$ and $S_{13}$ are algebraic.

The $4249$-dimensional simple module $S_{12}$ has an $8$-dimensional source $A_{12}$, and there are exactly $27$ elements of $\mc T(A_{12})$, the largest of which has dimension $80$, so that $S_{12}$ is algebraic. This completes the proof.
\end{pf}

The remaining simple modules for $p=5$ are $S_7$ and $S_8$, of dimension $925$, $S_9$ and $S_{10}$, of dimension $3197$, and $S_{14}$, of dimension $10860$. The modules $S_7$, $S_8$ and $S_{14}$ have dimension a multiple of $5$, and are non-periodic, so are likely to be non-algebraic, in accordance with Conjecture \ref{conj:periffalg}: the sources of $S_7$ and $S_8$ have dimension $75$, and the source of $S_{14}$ has dimension $110$, and at the moment we cannot prove whether these are algebraic or not.

The modules $S_9$ and $S_{10}$ are more complicated: they share a $47$-dimensional source $A_9$, and we have
\[ S^2(A_9)=k\oplus X_1\oplus A_4\oplus B_1\oplus 38\cdot\proj(k)\text{ and }\Lambda^2(A_9)=X_3\oplus A_5\oplus A_9\oplus A_{11}\oplus 34\cdot\proj(k).\]
Here, $X_i$ is the sum over all (six) subgroups $Q$ of $P$ of order $5$ of the $i$-dimensional indecomposable module for $Q$, induced to $P$, and $B_1$ is a $118$-dimensional simple module. Hence by Lemma \ref{lem:Mplusalgebraic}, if $B_1$ is algebraic then $A_9$ is algebraic. Taking the tensor square of $B_1$ is more difficult, and we just do the symmetric square. This decomposes as
\[ S^2(B_1)=k\oplus 2\cdot X_1\oplus X_3\oplus A_4\oplus B_1\oplus B_2\oplus B_3\oplus B_4\oplus 256\cdot\proj(k).\]
In this decomposition: $B_2$ is a $61$-dimensional indecomposable module, which appears as a summand in $S^2(A_7)$, hence algebraic; $B_3$ is a sum of four $10$-dimensional periodic modules, each of which is easily proved to be algebraic; and $B_4$ is a $222$-dimensional indecomposable module. The module $B_4$ is much more difficult to analyze, and at present we cannot get any information about this module. However, with the plethora of algebraic modules appearing in its tensor powers, the author believes that $S_9$ and $S_{10}$ are indeed algebraic.

\bibliography{references}

\providecommand{\bysame}{\leavevmode\hbox to3em{\hrulefill}\thinspace}
\providecommand{\MR}{\relax\ifhmode\unskip\space\fi MR }
\providecommand{\MRhref}[2]{%
  \href{http://www.ams.org/mathscinet-getitem?mr=#1}{#2}
}
\providecommand{\href}[2]{#2}
\begin{thebibliography}{10}

\bibitem{alperin1976b}
Jonathan Alperin, \emph{On modules for the linear fractional groups},
  International Symposium on the Theory of Finite Groups, 1974, Tokyo (1976),
  157--163.

\bibitem{alperin1979}
\bysame, \emph{Projective modules for {$\SL(2,2^n)$}}, J.~Pure Appl.~Algebra
  \textbf{15} (1979), 219--234.

\bibitem{archer2008}
Louise Archer, \emph{On certain quotients of the {G}reen rings of dihedral
  $2$-groups}, J.~Pure Appl.\ Algebra \textbf{212} (2008), 1888--1897.

\bibitem{auslandercarlson1986}
Maurice Auslander and Jon Carlson, \emph{Almost-split sequences and group
  rings}, J.~Algebra \textbf{103} (1986), 122--140.

\bibitem{bensoncarlson1986}
David Benson and Jon Carlson, \emph{Nilpotent elements in the {G}reen ring},
  J.~Algebra \textbf{104} (1986), 329--350.

\bibitem{conlon1966}
Samuel~B. Conlon, \emph{The modular representation algebra of groups with
  {S}ylow $2$-subgroup ${Z}_2\times {Z}_2$}, J.~Austral.~Math.~Soc. \textbf{6}
  (1966), 76--88.

\bibitem{craven2010un3}
David~A.\ Craven, \emph{Sources of simple modules for weight $2$ blocks of
  symmetric groups}, Submitted.

\bibitem{craventhesis}
\bysame, \emph{Algebraic modules for finite groups}, Ph.D. thesis, University
  of Oxford, 2008.

\bibitem{craven2009}
\bysame, \emph{Simple modules for groups with abelian {S}ylow $2$-subgroups are
  algebraic}, J.\ Algebra \textbf{321} (2009), 1473--1479.

\bibitem{craven2011}
\bysame, \emph{Algebraic modules and the {A}uslander--{R}eiten quiver}, J.\
  Pure Appl.\ Algebra \textbf{215} (2011), 221--231.

\bibitem{cekl2008}
David~A.\ Craven, Charles~W.\ Eaton, Radha Kessar, and Markus Linckelmann,
  \emph{The structure of blocks with {K}lein four defect group}, Math.\ Z., to
  appear.

\bibitem{donkin1993}
Stephen Donkin, \emph{On tilting modules for algebraic groups}, Math. Z.
  \textbf{212} (1993), 39--60.

\bibitem{dotyhenke2005}
Stephen Doty and Anne Henke, \emph{Decomposition of tensor products of modular
  irreducibles for {${\SL}_2$}}, Quart.\ J.\ Math. \textbf{56} (2005),
  189--207.

\bibitem{erdmann1977}
Karin Erdmann, \emph{Principal blocks of groups with dihedral {S}ylow
  $2$-subgroups}, Comm.~Algebra \textbf{5} (1977), 665--694.

\bibitem{erdmann1979}
\bysame, \emph{On $2$-blocks with semidihedral defect groups}, Trans.\ Amer.\
  Math.\ Soc. \textbf{256} (1979), 267--287.

\bibitem{erdmannhenke2002}
Karin Erdmann and Anne Henke, \emph{On {R}ingel duality for {S}chur algebras},
  Math. Proc. Camb. Phil. Soc. \textbf{132} (2002), 97--116.

\bibitem{feit1980}
Walter Feit, \emph{Irreducible modules of $p$-solvable groups},
  Proc.~Sympos.~Pure Math (Santa Cruz, 1979) \textbf{37} (1980), 405--412.

\bibitem{feit}
\bysame, \emph{The representation theory of finite groups}, North--Holland,
  Amsterdam--New York, 1982.

\bibitem{fongharris1993}
Paul Fong and Morton Harris, \emph{On perfect isometries and isotypies in
  finite groups}, Invent. Math. \textbf{114} (1993), 139--191.

\bibitem{abc}
Christoph Jansen, Klaus Lux, Richard Parker, and Robert Wilson, \emph{An atlas
  of {B}rauer characters}, Oxford University Press, New York, 1995.

\bibitem{jennings1941}
Stephen Jennings, \emph{The structure of the group ring of a $p$-group over a
  modular field}, Trans.\ Amer.\ Math.\ Soc. \textbf{50} (1941), 175--185.

\bibitem{kawatamichleruno2001}
Shigeto Kawata, Gerhard Michler, and Katsuhiro Uno, \emph{On
  {A}uslander-{R}eiten components and simple modules for finite groups of {L}ie
  type}, Osaka J.~Math. \textbf{38} (2001), 21--26.

\bibitem{kleidman1988}
Peter Kleidman, \emph{The maximal subgroups of the {C}hevalley groups
  {$G_2(q)$} with {$q$} odd, the {R}ee groups {$^2G_2(q)$}, and their
  automorphism groups}, J.\ Algebra \textbf{117} (1988), 30--71.

\bibitem{koshitanikunugi2002}
Shigeo Koshitani and Naoko Kunugi, \emph{Brou\'e's conjecture holds for
  principal $3$-blocks with elementary abelian defect group of order 9}, J.\
  Algebra \textbf{248} (2002), 575--604.

\bibitem{koshitanikunugiwaki2004}
Shigeo Koshitani, Naoko Kunugi, and Katsushi Waki, \emph{Brou\'e's abelian
  defect group conjecture holds for the {H}eld group and the sporadic {S}uzuki
  group}, J.\ Algebra \textbf{279} (2004), 638--666.

\bibitem{koshitanimiyachi2000}
Shigeo Koshitani and Hyoue Miyachi, \emph{The principal $3$-blocks of four- and
  five-dimensional projective special linear groups in non-defining
  characteristic}, J.\ Algebra \textbf{226} (2000), 788--806.

\bibitem{kovacs1981un}
L{\'a}szl{\'o} Kov{\'a}cs, \emph{Some indecomposables for $\mathrm{SL}_2$},
  Research report no. 11, Mathematics Research Report Series of the Australian
  National University, 1981.

\bibitem{liebeckseitz1996}
Martin~W. Liebeck and Gary~M. Seitz, \emph{Reductive subgroups of exceptional
  algebraic groups}, Mem.\ Amer.\ Math.\ Soc. \textbf{121} (1996), no.~580,
  vi+111.

\bibitem{puig1990}
Llu\'is Puig, \emph{Alg\`ebres de source de certains blocs des groupes de
  {C}hevalley}, Ast\'erisque \textbf{181-182} (1990), 221--236.

\bibitem{wilsonrob}
Robert~A.\ Wilson, \emph{The finite simple groups}, Graduate Texts in
  Mathematics, vol. 251, Springer-Verlag London Ltd., London, 2009.

\end{thebibliography}

\end{document}